\documentclass[12pt,leqno]{article}
\usepackage[dvips]{graphicx}
\usepackage{amsmath,amssymb,amsthm,amscd}
\usepackage[mathscr]{eucal}
\usepackage{amsfonts}
\begin{document}
\parskip=6pt
\newtheorem{prop}{Proposition}
\numberwithin{prop}{section}
\newtheorem{thm}{Theorem}
\numberwithin{thm}{section}
\newtheorem{corr}{Corollary}
\numberwithin{corr}{section}
\newtheorem{lemma}{Lemma}
\newtheorem{rmk}{Remark}
\newtheorem{defn}{Definition}
\numberwithin{defn}{section}
\numberwithin{lemma}{section}
\numberwithin{equation}{section}
\newcommand{\ds}{\displaystyle}
\newcommand{\va}{\vartheta}
\newcommand{\ep}{\epsilon}
\newcommand{\bC}{\mathbb C}
\newcommand{\bK}{\mathbb K}
\newcommand{\bP}{\mathbb P}
\newcommand{\bN}{\mathbb N}
\newcommand{\bA}{\mathbb A}
\newcommand{\bR}{\mathbb R}
\newcommand{\cH}{\mathcal H}
\newcommand{\pa}{\partial}
\newcommand{\ue}{\underline \ell}
\newcommand{\uh}{\underline h}
\newcommand{\uk}{\underline k}
\newcommand{\uu}{\underline u}
\newcommand{\ux}{\underline x}
\newcommand{\uy}{\underline y}
\newcommand{\uz}{\underline z}
\newcommand{\up}{\underline \pi}
\newcommand{\uP}{\underline P}
\newcommand{\urh}{\underline\rho}
\newcommand{\wta}{\widetilde A}
\newcommand{\wtp}{\widetilde p}
\newcommand{\wtk}{\widetilde {\underline{k}}}
\newcommand{\rp}{\operatorname{prim}}
\newcommand{\topdim}{\operatorname{top\;dim}}
\newcommand{\rk}{\operatorname{ker}}
\newcommand{\rn}{\operatorname{hull}}
\newcommand{\rd}{\operatorname{dim}}
\newcommand{\rt}{\operatorname{tsr}}
\newcommand{\rR}{\operatorname{RR}}
\newcommand{\rra}{\operatorname{rank}}
\newcommand{\rma}{\operatorname{max}}
\newcommand{\rmi}{\operatorname{min}}
\newcommand{\rdim}{\operatorname{dim}}
\newcommand{\la}{\lambda}
\newcommand{\bon}{{\bf{1}}}
\newcommand{\rht}{\rho_t(I)}
\newcommand{\opi}{\overline\pi}

\renewcommand\qed{ }
\begin{titlepage}
\title{\bf On Higher Real and Stable Ranks for $CCR$ $C^*-$algebras}
\author{Lawrence G. Brown}
\end{titlepage}
\date{}
\maketitle
\abstract
We calculate the real rank and stable rank of $CCR$ algebras which either have only finite dimensional irreducible representations or have finite topological dimension. We show that either rank of $A$ is determined in a good way by the ranks of an ideal $I$ and the quotient $A/I$ in four cases: When $A$ is $CCR$; when $I$ has only finite dimensional irreducible representations; when $I$ is separable, of generalized continuous trace and finite topological dimension, and all irreducible representations of $I$ are infinite dimensional; or when $I$ is separable, stable, has an approximate identity consisting of projections, and has the
corona factorization property. We also present a counterexample on higher ranks of $M(A)$, $A$ subhomogeneous, and a theorem of P. Green on generalized continuous trace algebras.

\endabstract



\section{Introduction.}
Rieffel [Ri] defined the (topological) stable rank, tsr$(A)$, of a $C^*$--algebra $A$, which by [HV] is the same as the Bass stable rank. Pedersen and I [BP1] defined the real rank, RR$(A)$, in an analogous way. A number of authors have given calculations of one or both of these ranks for naturally arising classes of $C^*$-algebras. The bibliography of [AK] contains a large list of such papers.  Many of these works have used theorems about rank for special classes of $CCR$ algebras or for extensions where the ideal is in a special class of $CCR$ algebras. This paper arises from [BP2, Theorem 2.12], which makes it possible to generalize some of these theorems. In particular, I am consciously generalizing results of Nistor [Ni2]. Some of the lemmas are stated for algebras which may not be $CCR$, or even type $I$, and it is possible that these lemmas, as well as the theorems, could be useful for calculating the ranks of additional naturally arising $C^*$--algebras. Although most of the results are stated in a way that includes the low ranks, stable rank one and real rank zero, the low rank cases were already known.

The reason for drawing lines between stable rank one and all higher values of stable rank and between real rank zero and all higher values of real rank, in the phrases ``low rank'' and ``higher rank'', is that the low ranks have different formal properties from the higher ranks. For example, the low ranks are invariant under Rieffel--Morita equivalence, whereas tsr$(A\otimes\bK)=2$ whenever tsr$(A)>1$ and RR$(A\otimes\bK)=1$ whenever RR$(A)>0$. Another example is found in the relation between rank and extensions, where our knowledge is far from complete.

Rieffel [Ri]  showed that if $I$ is a (closed, two--sided) ideal of a $C^*$--algebra $A$, then either
\item{(1)} tsr$(A)=\rma(\rt(I), \rt(A/I))$, or
\item{(2)} tsr$(A)=1 +\rma(\rt(I), \rt(A/I))$.

\noindent In the case tsr$(I)=\rt(A/I)=1$, (1) holds if and only if a natural lifting condition is satisfied. And this lifting condition is equivalent to the vanishing of the boundary map, $\pa_1:K_1(A/I)\to K_0(I)$. This $K$-theoretic criterion was first obtained by G. Nagy, cf. [Ni1, Lemma 3], and an alternate proof was published in [Na1, Corollary 2]. But for the higher rank case, so far as I know, no liftability criterion for (1) has been found except in special cases, and also no example has been found where max(tsr$(I),{\rm{tsr}}(A/I))>1$ and (2) holds.

If RR$(I)={\rm{RR}}(A/I)=0$, then RR$(A)=0$ if and only if a natural lifting condition is satisfied. And this lifting condition is equivalent to the vanishing of the other boundary map, $\pa_0:K_0(A/I)\to K_1(A)$. The $K$--theoretic criterion was first proved by S. Zhang, cf. [BP1, Propositions 3.14 and 3.15]. But we know much less about the higher real rank case in general. It is obvious that RR$(A/I)\leq {\rm{RR}}(A)$, and N. Hassan, [H, Theorem 1.4], showed that RR$(I)\leq {\rm{RR}}(A)$.

The two calculations mentioned in the abstract of the ranks of $CCR$ algebras are Theorems 3.9 and 3.10 below. The four results on ranks of extensions are Theorems 3.6, 3.11 3.12, and Corollary 3.15, and Corollary 3.14 illustrates the use of bootstrap methods to get formally stronger results.  Corollary 3.15 is an
afterthought which makes no use of $CCR$ algebras or [BP2, Theorem 2.12]. 
The counterexample on the ranks of $M(A)$, along with some related remarks and questions, is given in 3.16.

In both 3.11 and 3.12, the ideal $I$ is of generalized continuous trace (GCT). In the one case $I$ has only finite dimensional irreducible representations, and in the other, only infinite dimensional. There isn't any obvious way to reduce the study of arbitrary GCT algebras to these two cases. Section 4 contains new characterizations of separable GCT algebras, all but one unpublished results of P. Green [G] included here with his permission. Green's main result is that a separable $C^*$--algebra is GCT if and only if it is stably isomorphic to a $C^*$--algebra with only finite dimensional irreducibles. Although these characterizations aren't needed for the main results, they provide an interesting context. Also, Green's work, helped me to develop the perspective needed for my work. Finally, if GCT turns out to be the ``right'' hypothesis, within the class of $CCR$ 
algebras, for results like 3.11 and 3.12, perhaps the material in Section 4 will be helpful in getting better results.

I also thank R. Archbold for helpful comments.

\noindent\textbf{Bibliographical and Personal Remark.}
\noindent
I obtained the stable rank versions, in the separable case, of 3.8 and 3.9, and probably also 3.10 and 3.11, when I was working with Gert Pedersen in the late 1990's. We were working on [BP2] among other things, and I meant for these results to go into [BP2]. But Gert didn't want the paper to include results on higher rank unless they followed either from the same proofs as our low rank results or with minimal additional effort. We therefore agreed that I would publish these results separately after [BP2] was complete. I then put this subject aside, apparently without making notes of the statements or proofs. When I returned to the subject in connection with the completion of [BP2], I obtained better results, in particular the real rank versions. The paper [BP2] is the second--to--last of Gert's and my joint papers. Working with Gert was one of the best experiences of my life.

\section{Preliminaries.}
{\bf 2.1. Definitions.} If $A$ is a unital $C^*$--algebra and $\underline x=(x_1,\dots, x_n)$ is in $A^n$, then $\ux$ is \textit{unimodular} if it is left invertible when considered as an $n\times 1$ matrix. It is equivalent to require that $\sum^n_1\,x^*_ix_i$ be invertible or that $\{x_1,\dots,x_n\}$ generate $A$ as a left ideal. Then tsr$(A)$ is the smallest $n$ such that unimolular $n$--tuples are dense in $A^n$ and RR$(A)$ is the smallest $n$ such that unimodular $(n+1)$--tuples $\ux$ for which each $x_i$ is self--adjoint, are dense in $(A_{\mathrm{sa}})^{n+1}$. If no such $n$ exists, then the rank of $A$ is $\infty$. Thus $1\leq{\rm{tsr}}(A)\leq\infty$ and $0\leq{\rm{RR}}(A)\leq\infty$. If $A$ is non--unital, then define tsr$(A)={\rm{tsr}}(\wta)$ and RR$(A)={\rm{RR}}(\wta)$, where $\wta$ is the unitization of $A$.

\medskip
\noindent
{\bf 2.2. The Primitive Ideal Space}.
The primitive ideal space of $A$ is denoted prim$(A)$. Even when $A$ is type $I$, so that prim$(A)$ is identified with the spectrum of $A$, I continue to use this notation. If $F$ is a closed subset of prim$(A)$, then ker$(F)$ is the ideal $I$ defined by $I=\cap_{P\in F}P$, and $F={\rm{hull}} (I)=\{P\in{\rm{prim}}(A):P\supset I\}$. Also prim$(A/I)$ is identified with $F$, and prim$(I)$ is identified with prim$(A)\setminus F$. If $S$ is a locally closed subset of prim$(A)$, i.e., $S=F\cap G$ with $F$ closed and $G$ open, then $S$ is identified with prim$(I/J)$, where $I$ and $J$ are ideals such that $I\supset J$ and $S={\rm{hull}}(J)\setminus{\rm{hull}}(I)$. Although $I$ and $J$ are not uniquely determined by $S$, the quotient $I / J$ is determined up to canonical isomorphism. Thus $I / J$ may be denoted by $A(S)$. In particular, for $F$ closed $A(F)=A/{\rm{ker}}(F)$, and for $G$ open $A(G)={\rm{ker}}(\rp(A)\setminus G)$. It follows from results stated above that $\rt(A(S))\leq\rt(A)$ and RR$(A(S))\leq{\rm{RR}}(A)$.

\medskip
\noindent
{\bf 2.3. The Countable Sum Theorem}.
Parts (i) and (ii) of [BP2, Theorem 2.12] can be stated as follows:

\medskip
\noindent
(CST)\quad If $\rp(A)=\bigcup^\infty_{n=1} F_n$, \quad{\rm{where\; each}}\; $F_n$\;{\rm{is\; closed,\; then

\qquad$\rt(A)={\rm{sup}}_n\{\rt(A(F_n))\}$\;and\;$\rR(A)={\rm{sup}}_n\{{\rm{RR}}(A(F_n))\}$.

\medskip
\noindent
{\bf 2.4. Definitions.} The concept of generalized continuous trace (GCT) was defined by Dixmier [D2, $\oint$10], cf. also [D3, 4.7.12]. Let $J(A)$ denote the closure of the set of continuous trace elements of $A$. Then $J(A)$ is the largest ideal of $A$ such that $J(A)$ has continuous trace as a $C^*$--algebra and every compact subset of $\rp(J(A))$ is closed in $\rp(A)$. (In general there is no largest continuous trace ideal.) The continuous trace composition series is $\{J_\alpha:0\leq\alpha\leq\beta\}$, where $\beta$ is an ordinal number, $J_0=0, J_\la=(\cup_{\alpha<\la} J_\alpha)^-$ for $\la$ a limit ordinal, $J_{\alpha+1}/J_\alpha=J(A/J_\alpha)\neq 0$ for $\alpha<\beta$, and $J(A/J_\beta)$=0. Then $A$ is GCT if and only if $J_\beta=A$. Although every type $I$ $C^*$--algebra has a composition series with continuous trace quotients, every GCT $C^*$--algebra is $CCR$. Dixmier proved that GCT algebras are distinguished from other type $I$ $C^*$--algebras by the topology of their spectra.

\noindent
{\bf 2.5. Topological Dimension}.
A topological space is called {\textit{almost}} {\textit{Hausdorff}} if every non--empty closed subset $F$ contains a non--empty relatively open subset which is Hausdorff in the relative topology. Thus $\rp(A)$ is almost Hausdorff whenever $A$ is type $I$. In [BP2] $\topdim(A)$ was defined for $C^*$--algebras $A$ with almost Hausdorff primitive ideal space as follows: $\topdim(A)=\sup_K\{\rdim K\}$, where $K$ ranges through compact Hausdorff (locally closed) subsets of $\rp(A)$ and dim denotes covering dimension. Thus $\topdim(A)$ is a topological property of $\rp(A)$, but it is not the same as dim$(\rp(A))$. If $\rp(A)$ is Hausdorff, then $\topdim(A)={\rm{loc\; dim}}(\rp(A))$, which is the same as dim$(\rp(A)\cup\,\{\infty\})$, where $\rp(A) \cup\{\infty\}$ is the one-point compactification, (and the same as dim$(\rp(A))$ if $\rp(A)$ is $\sigma$--compact). It was shown in [BP2] that $\topdim(A)$ behaves well under extensions and composition series, and it was explained why $\topdim(A)$ is a better choice than $\rdim(\rp(A))$ when they differ.

The following easy lemma will be used in the proof of the real rank case of 3.12.

\medskip
\noindent   
{\bf Lemma\ 2.6.\ }\textit{Let $A$ be a non--zero unital $C^*$--algebra, and let
$\uh$ be an $n-$tuple in $(A_{sa})^n$, where $n\geq 2$. Then the $n\times n$ matrix $(h_ih_j)$ is not invertible.}

\begin{proof} Regard $\uh$ as an $n\times 1$ matrix, so that the matrix in question is $\uh\,\uh^*$. If $A$ is a unital subalgebra of $B(\cH)$, then $\uh$ may 
be regarded as an operator from $\cH$ to $\cH\oplus\cdots\oplus\cH$.  Then if $\uh\,\uh^*$ is invertible, $\uh$ must be surjective.  It follows that each $h_i$ is 
surjective and (since $n>1$) no $h_i$ is injective.  This is absurd, since
$h_i$ is self--adjoint.
\end{proof}

\section{Main Results.}
Many of the proofs are essentially the same for the stable rank and real rank cases. The notation rank$(A)$ will be used to denote either tsr$(A)$ or RR$(A)$ in such proofs.

\begin{lemma}\textit{ Let $I$ be an ideal of a $C^*$--algebra $A$. Assume that prim$(I)$ is Hausdorff and each compact subset of prim$(I)$ is closed in prim$(A)$. Then $\rt(A)={\rm{max}}(\rt(I),\rt(A/I))$, and RR$(A)={\rm{max}}({\rm{RR}}(I),{\rm{RR}}(A/I))$.}
\end{lemma}

\begin{proof} This can be deduced from [Sh, Proposition 3.15] and 
its real rank counterpart, 
[O1, Lemma 1.9].  Let $\Lambda$ be the set of relatively compact open 
subsets of $\rp(I)$, and let $J_{\lambda}$ be the corresponding ideal for 
$\lambda\in\Lambda$.  (Thus, in the notation of 2.2, $J_{\lambda}=I(\lambda)=A(\lambda)$.)
Then $\{J_{\lambda}\}$ is upward directed, $I=(\bigcup_{\lambda}J_{\lambda})^-$, 
and the results cited tell us that $\rra(A)$ is the larger of $\rra(A/I)$ and 
$\mathrm{sup}_{\lambda}\{\rra(A/J_{\lambda}^{\perp})\}$.  But $A/J_{\lambda}^{\perp}=A(\overline\lambda)$, which by hypothesis is a quotient of $I$.
\end{proof}

\medskip
\noindent   
{\bf Remark.} If $I$ is $\sigma$--unital, the Lemma can also be deduced from 
(CST), since then $\rp(I)$ is an $F_{\sigma}$ in $\rp(A)$. As noted in 
[BP2, Remark 3.11], Sheu's Technical Proposition,  [Sh, Proposition 3.15], helped to inspire (CST) and in turn could be deduced from (CST). 

\medskip
\noindent
{\bf Definition 3.2.} If $X$ is a primitive ideal space, then an \textit{$FD$--like decomposition} of $X$ is a family $\{H_1, H_2,\dots\}$ of locally closed subsets of $X$ such that:
\item{(i)} $X=\bigcup_nH_n,\, H_n\cap H_m=\emptyset$ if $n\neq m$.
\item{(ii)} Each $H_n$ is Hausdorff.
\item{(iii)} Every compact subset of $H_n$ is closed in $X$.
\item{(iv)} $F_n=\bigcup_{k=1}^n H_k$ is closed.

\noindent
The terminology is explained by the following result, which is stated only for reference, since it is well known. 
\medskip

\noindent
{\bf Proposition\ 3.3.\ }\textit{Let $A$ be a $C^*$--algebra all of whose irreducible representations are finite dimensional, and let $H_n=\{P\in\rp(A):A/P\cong \mathbb M_n\}$. Then $\{H_n\}$ is an $FD$--like decomposition of $\rp(A)$.}

\medskip
\noindent
{\bf Lemma\ 3.4.\ }\textit{If $\{H_n\}$ is an $FD$--like decomposition of $\rp(A)$, then $\rt(A)=\sup_n\{\rt(A(H_n))\}$ and $\rR(A)=\sup_n\{\rR(A(H_n))\}$.}

\begin{proof}That rank$(A(H_n))\leq\rra(A)$ is clear.  By (CST) it is sufficient to show $\rra(A(F_n))\leq \sup_m\{\rra(A(H_m))\}$, $\forall n$. This is done by induction on $n$, the case $n=1$ being obvious. For $n>1$, $A(F_n)$ contains $A(H_n)$ as an ideal, and the quotient is $A(F_{n-1})$. Thus the result follows from Lemma 3.1.
\end{proof}

\medskip
\noindent
{\bf Lemma\ 3.5.\ }\textit{If $I$ is an ideal of $A$ and if $\rp(A)$ has an $FD$--like demomposition, then $\rt(A)=\rma(\rt(I),\rt(A/I))$, and $\rR(A)=\rma(\rR(I), \rR(A/I)$.}

\begin{proof}Let the $FD$--like decomposition be $\{H_n\}$. Then by Lemma 3.4 it is enough to show that for each $n$ we have 
$\rra(A(H_n))=\rma(\rra(A(H_n\cap \rp(I))), \rra(A(H_n\cap \rn(I))))$. But this follows directly from Lemma 3.1, since $A(H_n)$ has a Hausdorff primitive ideal space.
\end{proof}

\medskip
\noindent
{\bf Theorem\ 3.6.\ }\textit{If $A$ is a $CCR$ $C^*-$algebra and $I$ a closed two--sided ideal, then}
\item{(i)\ }  $\rt(A)=\rma(\rt(I),\rt(A/I))$, \textit{and}
\item{(ii)\ } $\rR(A)=\rma(\rR(I),\rR(A/I))$.

\begin{proof}We can write $A=(\bigcup B_i)^-$, where $\{B_i\}$ is an upward directed family of hereditary $C^*$--subalgebras, and each $B_i$ has only finite dimensional irreducible representations. This can be deduced from the theory of the Pedersen ideal, $K(A)$, which is the minimum dense two--sided ideal of $A$. Each $B_i$ will be the hereditary $C^*$--subalgebra generated by a finite subset of $K(A)$. Since $\pi(x)$ has finite rank for each irreducible $\pi$ and each $x$ in $K(A)$, $B_i$ has the required property. Let $J_i$ be the ideal generated by $B_i$. Because of the compatibility of rank with direct limits, it is enough to show rank$(J_i)=\rma(\rra(J_i\cap I)$, $\rra(J_i/J_i\cap I))$ for each $i$. Since prim$(J_i)$ is homeomorphic to prim$(B_i)$, this follows from Proposition 3.3 and Lemma 3.5.
\end{proof}

\medskip
\noindent
{\bf Corollary\ 3.7.\ }\textit{If $A$ is a $CCR$ $C^*$--algebra and $\{I_\alpha: 0\leq\alpha\leq \beta\}$ is a composition series for $A$, then}
\item{(i)\ }  $\rt(A)=\underset{\alpha<\beta}{\sup}\{\rt(I_{\alpha +1}/I_\alpha)\}$, \textit{and}
\item{(ii)\ }  $\rR(A)=\underset{\alpha<\beta}{\sup}\{\rR(I_{\alpha+1}/I_\alpha)\}$.

\begin{proof}Let $t$ be the sup. We prove by transfinite induction that rank$(I_\alpha)\leq t$. If $\alpha$ is a limit ordinal and the result is true for $\gamma<\alpha$, then it is true for $\alpha$ by a direct limit argument. And if $\alpha=\gamma+1$ and the result is true for $\gamma$, then the theorem implies it for $\alpha$.
\end{proof}

\noindent
{\bf Proposition\ 3.8.\ }\textit{If $A$ is $n$--homogeneous and top dim$(A)=d$, then $\rt(A)=\lceil{2n-1+d\over 2n}\rceil=\lfloor{4n-2+d\over 2n}\rfloor$, and $\rR(A)=\lceil{d\over 2n-1}\rceil=\lfloor{d+2n-2\over 2n-1}\rfloor$.}

\begin{proof}Let $\rp(A)=X$, so that $A$ is the algebra of continuous sections vanishing at $\infty$ of a locally trivial $\mathbb M_n$--bundle over $X$. Note that the formula is known if $X$ is compact and $A=C(X)\otimes M_n$ by [Ri] for the stable rank case and [BEv] for the real rank case. Each compact subset $K$ of $X$ can be written $K=F_1\cup\dots\cup F_k$ where each $F_i$ is closed and the bundle is trivial over $F_i$. Since dim$(K)={\rm{\max}}^k_{i=1}(\rdim(F_i))$ and rank$(A(K))={\rm{\max}}^k_{i=1}\rra(A(F_i))$ (by (CST) or [Sh] and [O1]), it is clear that rank$(A)$ is at least the number given. For the reverse inequality write $X=\bigcup U_i$ where $\{U_i\}$ is an upward directed family of $\sigma$--compact open subsets. Then by a direct limit augument, rank$(A)\leq\sup_i\{\rra(A(U_i))\}$. Each $U_i$ is a countable union of compact subsets on which the bundle is trivial. Thus the result follows from (CST).
\end{proof}

\medskip
\noindent
{\bf Theorem\ 3.9.\ }\textit{Let $A$ be a $C^*$--algebra with only finite dimensional irreducible representations and $H_n=\{P\in\rp(A):A/P\cong \mathbb M_n\}$. Then if $\topdim(A(H_n))\mathrm{ (=loc\; dim}(H_n))=d_n\; (d_n=0\;\mathrm{if}\;H_n=\emptyset)$, we have}
\item{(i)\ }  $\rt(A)=\sup_n\{\lceil{2n-1+d_n\over 2n}\rceil\}$, \textit{and}
\item{(ii)\ } $\rR(A)=\sup_n\{\lceil{d_n\over 2n-1}\rceil\}$.

\begin{proof}Combine 3.3, 3.4, and 3.8.
\end{proof}

\medskip
\noindent
{\bf Theorem\ 3.10.\ }\textit{Let $A$ be a $CCR$ $C^*-$algebra, and suppose
that $d=\topdim(A)<\infty$. Let $H_n=\{P\in\rp(A):A/P\cong \mathbb M_n\}$, and
let $d_n=\topdim(A(H_n)) (=\rm{loc}\;\rdim(H_n))$.}
\item{(i)\ }   \textit{If} $d\leq 1$, \textit{then} $\rt(A)=1$.
\item{(ii)\ }  \textit{If} $d>1$, \textit{then} $\rt(A)=\sup_n\{\rma(\lceil{2n-1+d_n\over 2n}\rceil, 2)\}$.
\item{(iii)\ } \textit{If} $d=0$, \textit{then} $\rR(A)=0$.
\item{(iv)\ }  \textit{If} $d>0$, \textit{then} $\rR(A)=\sup_n\{\rma(\lceil{d_n\over 2n-1}\rceil, 1)\}$.

\begin{proof}It is already known that $\rt(A)=1$ if and only if $d\leq 1$ and $\rR(A)=0$ if and only if $d=0$. For the stable rank case this is [BP2, Theorem 5.6]. For the real rank case, it follows from [BP2, Proposition 5.1], but, as explained in [BP2], was previously known from Bratteli and Elliott [BEl] if $A$ is separable. This proves parts (i) and (iii) as well as the fact that rank$(A)$ is at least the number indicated in parts (ii) and (iv).

Let $N$ be a positive integer such that ${2N-1+d\over 2N}\leq 2$ and ${d\over 2N-1}\leq 1$, let $F_{N-1}\subset\rp(A)$ be defined as above, and let $I=\rm{ker}(F_{N-1})$. Then $\rra(A)=\rma(\rra(I),\;\rra(A/I))$ by 3.6, and $\rra(A/I)$ can be computed by 3.9, since $A/I$ is subhomogeneous (prim$(A/I)=F_{N-1})$. Thus all irreducible representations of $I$ have dimension at least $N$, and it is sufficient to show $\rt(I)\leq 2$ and $\rR(I)\leq 1$. It can be shown that $I=(\bigcup_iB_i)^-$, where $\{B_i\}$ is an upward directed family of hereditary $C^*$--subalgebras each of whose irreducible representations has finite dimension at least $N$. But a slightly roundabout approach seems less technical.

As in the proof of 3.6, write $I=(\bigcup J_i)^-$, where $\{J_i\}$ is an upward directed family of ideals such that each prim$(J_i)$ has an $FD$--like decomposition. Since it is sufficient to show tsr$(J_i)\leq 2$ and $\rR(J_i)\leq 1$, we may assume prim$(I)$ has an $FD$--like decomposition. Then using a decomposition and Lemma 3.4, we reduce to the case where prim$(I)$ is Hausdorff. If $X=\rp(I)$, another direct limit augument reduces to the case where $X$ is $\sigma$--compact, and then an application of (CST) reduces to the case where $X$ is compact.

So after this final reduction we have a new $CCR\; C^*$--algebra, $A_1$, such that top dim$(A_1)\leq d$, all irreducible representations of $A_1$ have dimension at least $N$, and prim$(A_1)$ is compact Hausdorff. Write $A_1=(\bigcup C_j)^-$, where $\{C_j\}$ is an upward directed family of hereditary $C^*$--subalgebras each of which has only finite dimensional irreducible representations. For each $j$ let $U_j=\{x\in\rp(A_1):\rdim\;{\pi_x}_{\,|C_j}\geq N\}$, where $\pi_x$ is an irreducible representation with kernel $x$. Then $\{U_j\}$ is an open cover of prim$(A_1)$. By compactness $U_{j_0}=\rp(A_1)$ for some $j_0$. Hence $j\geq j_0$ implies all irreducible representations of $C_j$ have dimension at least $N$, which implies by Theorem 3.9 that tsr$(C_j)\leq 2$ and $\rR(C_j)\leq 1$.
\end{proof}

The proof given for the next theorem is a slightly simplified version, suggested
by R. Archbold, of the original proof.

\medskip
\noindent
{\bf Theorem\ 3.11.\ }\textit{If $I$ is an ideal of the $C^*$--algebra $A$ such that all irreducible representations of $I$ are finite dimensional, then}
\item{(i)\ }  $\rt(A)=\rma(\rt(I),\rt(A/I))$, \textit{and}
\item{(ii)\ } $\rR(A)=\rma(\rR(I),\rR(A/I))$.

\begin{proof}Let $F_n=\{P\in\rp(A):\rd(A/P)\leq n^2\}$ for $n\geq 1$ and $F_0=\rn(I)$. Apply (CST) to $\rp(A)=\bigcup^\infty_{n=0} F_n$. Thus it is sufficient to prove $\rra(A(F_n))\leq\rma(\rra(I),\;\rra(A/I))$ for $n>0$. Since $A(F_n)$ is subhomogeneous, it follows from either Lemma 3.5 or Theorem 3.6 that $\rra(A(F_n))\leq\rma(\rra(I), \rra(A/I))$, since $A(F_n)$ has an ideal $J$ which
 is a quotient of $I$ such that $A(F_n)/J$ is a quotient of $A/I$.
\end{proof}

The statement of the next theorem does not include the known facts in the case tsr$(I)=\rt(A/I)=1$, which are instead reviewed in Remark 3.13 (ii). The statement does fully cover the case $\rR(I)=\rR(A/I)=0$, but the proof does not deal with this case. Instead a stronger result is proved in Remark 3.13 (iii). Most of the content of Remark 3.13 (iii) resides in the already known results cited there.

\medskip
\noindent
{\bf Theorem\ 3.12.\ }\textit{Let $I$ be an ideal of the $C^*$--algebra $A$ such that $I$ is separable, $I$ has generalized continuous trace, $\topdim(I)<\infty$, and all irreducible representations of $I$ are infinite dimensional. Then}
\item{(i)\ }$\rt(A)\leq\rma(2,\rt(A/I))$, \textit{and}
\item{(ii)\ }$\rR(A)=\rma(\rR(I),\rR(A/I))$.

\begin{proof}Let $\{J_\alpha:0\leq\alpha\leq\beta\}$ be the continuous trace composition series for $I$ defined in 2.4. Here $\beta$ is a countable ordinal number and $J_\beta=I$. Since each $J_{\alpha+1}/J_\alpha$ is a separable continuous trace $C^*$--algebra, then prim$(J_{\alpha+1}/J_\alpha)=\bigcup^\infty_{n=1}K_{\alpha,n}$, where the $K_{\alpha,n}$ are compact subsets such on each $K_{\alpha,n}$, $J_{\alpha+1}/J_\alpha$ is derived from a continuous field of Hilbert spaces. It then follows from the hypotheses and a result of Dixmier and Douady, [DD, Th\'{e}or\`{e}me 5], that each of these continuous fields is trivial. Moreover, each $K_{\alpha,n}$ is closed in prim$(I)$. Thus, after re--numbering, prim$(I)=\bigcup^\infty_{n=1}\,K_n$ where each $K_n$ is closed and compact Hausdorff and $I(K_n)\cong C(K_n)\otimes \bK$. Now let $F_n=K_n\cup\, \rn(I)\subset\rp(A)$, and apply (CST) to $\rp(A)=\bigcup^\infty_{n=1} F_n$. Thus we are reduced to the case $I=C(T)\otimes\bK$, where $T$ is compact, metrizable, and finite dimensional, Part (i) now follows directly from Nistor's result, [Ni2, Lemma 2].

For part (ii) we assume, as we may, that $A$ is unital and that $n=\rma(2, 1+\rR(A/I))<\infty$. Then we need to approximate a given tuple $\ux$ in $(A_{\mathrm{sa}})^n$ with a unimodular tuple in $(A_{\mathrm{sa}})^n$. If $A\subset B(\cH)$, then tuples will be regarded as operators from $\cH$ to $\cH^n=\cH\oplus\cdots\oplus\cH$. (The Hilbert space $\cH$ may be non--separable.) Let $\pi:A\to A/I$ be the quotient map and $\rho:A\to M(I)$ the natural map. Symbols such as $\up,\urh$ (respectively, $\pi_n,\rho_n)$ will denote the natural extensions to $A^n$ (respectively, $\mathbb M_n(A)$). Since $M(I)$ can be identified with the algebra of double--strongly continuous functions from $T$ to $B(\ell^2)$, the symbol $\rho_t(a)$, for example, will denote the value of $\rho(a)$ at the point $t$ in $T$.

If $\ep>0$, by the assumption on $\rR(A/I)$ there is a unimodular tuple $\up(\uy)$ with entries in $(A/I)_{\mathrm{sa}}$ such that $\|\up(\ux)-\up(\uy)\|<{\epsilon\over 2}$. By the properties of quotient norms we may assume $\|\ux-\uy\|<{\epsilon\over 2}$. Because $T$ is compact, $I$ has an approximate identity $(p_m)$ consisting of full projections. We claim that $\uy(\bon-p_m)$ is left invertible as an operator on $(\bon-p_m)\cH$ for $m$ large enough. It is sufficient to work with $|\uy|=(\sum^n_1\,y^*_iy_i)^{1\over 2}$. If $\delta>0$ is such that $\pi(|\uy|)\geq \delta\cdot\bon$, then $(|\uy|^2-\delta^2\cdot\bon)_-\,\in I$. Choose $m$ so that $\|(\bon-p_m)(|y|^2-\delta^2\cdot\bon)_-(\bon-p_m)\|<{3\over 4}\delta^2$. Then since $|\uy(\bon-p_m)|^2= (\bon-p_m)|\uy|^2(\bon-p_m)$, we conclude that $|\uy(\bon-p_m)|\geq{\delta\over 2}(\bon-p_m)$.

Let $p=p_m$ for $m$ as above and let $q$ in $\mathbb M_n(A)$ be the range projection of $\uy(\bon-p)$. We claim that $(\rho_n)_t(\bon_n-q)$ has infinite rank for each $t$ in $T$. This follows from Lemma 2.6 applied in $\rho_t(A)/\rho_t(I)$. Since $\rho_t(I)=\bK$ and $\rho_t(A)$ contains the identity of $B(\ell^2)$, this quotient is non--zero. Also note that $\urh_t(\uy(\bon-p))+\urh_t(I^n)=\urh_t(\uy)+\urh_t(I^n)$, since $p\in I$. Hence all entries are self--adjoint. Thus 
$(\rho_n)_t(\bon_n-q)\not\subset \mathbb M_n(\bK).$
Now results of Dixmier and Douady, [DD, Th\'{e}or\`{e}me 5 and Corollaire 3], imply that $\rho_n(\bon_n-q)=\sum^\infty_1\,r_m$ where the $r_m$'s are mutually orthogonal projections, each of which is Murray--von Neumann equivalent to $p$, and convergence is in the strict topology of $\mathbb M_n(M(I))=M(\mathbb M_n(I))$.

Operators from $\cH$ to $\cH^n$ will be represented as $2\times 2$ matrices relative to $\cH=(\bon-p)\cH\oplus p\cH$ and $\cH^n=q\cH^n\oplus (\bon_n-q)\cH^n$. If $\uz=\begin{pmatrix} a & b\\ c & d\end{pmatrix}$, and if $a$ is invertible (as an operator from $(\bon-p)\cH$ to $q\cH^n$), then, as is well known, $\uz$ is left invertible if and only if $d-ca^{-1}b$ is left invertible. If $\up(\uz)$ is unimodular, it is sufficient that $\urh(d-ca^{-1}b)$ be left invertible (since $\rk(\pi)\cap\rk(\rho)=0$), and for this it is sufficient that $r_m\urh(d-ca^{-1}b)$ be left invertible for one value of $m$. Of course, by construction $\uy=\begin{pmatrix}a_0 & b_0\\ 0 & d_0\end{pmatrix}$, where $a_0$ is invertible.

We will find a unimodular tuple $\uz$ such that $\|\uz-\uy\|<{\ep\over 2}$ and $\uz-\uy\in I^n$. We first choose an appropriate tuple $\uk=(k_1,\dots, k_n)$ in $(I^n)p$ and then take $\uz=\uy+\uk+\wtk$, where $\wtk=(k^*_1,\dots,k^*_n)$. Since $\|\wtk\|\leq n\|\uk\|$, one condition will be that $\|\uk\|<\ep/2(n+1)$. Let $\eta=\rmi(\ep/4(n+1), 1/2n, 1/4n\|a_0^{-1}\|)$. Let $\wtp=\mathrm{diag}\,(p,p,\dots,p)\in \mathbb M_n(I)$. Since $(r_m)$ converges strictly to 0, there is a value of $m$ such that $\|r_m\,d_0\|<\eta$ and $r_m\,\|\wtp\|<1/n(2+4\|a_0^{-1}\|(\|b_0\|+1))$. Then choose $\uu$ in $I^n$, such that $\uu^*\uu=p$ and $\uu\,\uu^*=r_m$, and let $\uk=\begin{pmatrix}0 & 0\\ 0 & \eta\uu-r_m d_0\end{pmatrix}$. (Note that $\rho$ is an isomorphism on $I$, so there is no need to distinguish $\uu$ from $\urh(\uu)$, $d_0$ from $\urh(d_0)$, or $r_m$ from $({\rho_n}_{|\mathbb M_n(I)})^{-1}r_m)$. Let $\wtk=\begin{pmatrix}a_1 & b_1\\ c_1 & d_1\end{pmatrix}$, so that 
$\|a_1\|, \|b_1\|, \|c_1\|, \|d_1\|\leq \|\wtk\|\leq n\|\uk\|<2n\eta$.
Then $\uz=\begin{pmatrix}a_0+a_1 &b_0+b_1\\ c_1 &d_0+d_1+\eta\uu-r_md_0\end{pmatrix}=\begin{pmatrix}a_2 &b_2\\ c_2 &d_2\end{pmatrix}$. Thus
\[ r_m\urh(d_2-c_2a_2^{-1}b_2)=\eta\uu+r_md_1-r_mc_1(\urh(a_0+a_1))^{-1}(b_0+b_1).\]
 Since $\uk=\uk p$, $\wtk=\wtp\wtk$, and hence 
\[\|r_m\wtk\|=\|r_m\wtp\wtk\|\leq\|r_m\wtp\|\;\|\wtk\|<2n\eta\|r_m\wtp\|.\]
 In particular $\|r_mc_1\|, \|r_md_1\| < 2n\eta\|r_m\wtp\|$. Also note that $\|a_1\|<1/2\|a_0^{-1}\|$, so that $a_0+a_1$ is invertible and $\|(a_0+a_1)^{-1}\|\leq 2 \|a_0^{-1}\|$. It is then routine to check that $\|r_md_1-r_mc_1(\urh(a_0+a_1))^{-1}(b_0+b_1)\|<\eta$. Hence $\uz$ is unimodular.
\end{proof}

\medskip
\noindent   
{\bf Remark\ 3.13.\ }(i) By Theorem 3.10 tsr$(I)=1$ or 2, according as 
$\topdim(I) \le 1$ or $\topdim(I)>1$, and $\rR(I)=0$ or 1, according as $\topdim(I)=0$ or $\topdim(I)>0$.

(ii) As previously mentioned, if tsr$(I)=\rt(A/I)=1$, then $\rt(A)$ can be determined using $K$--theory. The special assumptions on $I$ do not eliminate the need to look at the $K$--theory.

(iii) If $I$ is an arbitrary type $I$ $C^*$--algebra of real rank zero, or more generally if $I$ is any $AF$--algebra, then $\rR(A)=\rR(A/I)$. In fact Proposition 3.4 of Osaka's survey article [O2], which is obtained by combining Busby's analysis of extensions [Bu] with a pullback result of Nagisa, Osaka, and Phillips, [NOP, Proposition 1.6], states that $\rR(A)\leq\rma(\rR(M(I)), \rR(A/I))$.  
(The case where the max is $0$ was independently proved in 
[BP2, Corollary 4.4].) And a result of H. Lin, [L, Corollary 3.7], implies that $\rR(M(I))=0$ if $I$ is separable and $AF$. The fact that every separable type $I$ $C^*$--algebra of real rank zero is $AF$ follows from a result of Bratteli and Elliott, [BEl, \S7]. Finally, the separability hypothesis on $I$ can be removed via standard techniques for reducing to the separable case, cf. the proof of [BP1, Theorem 3.8]. Either the type $I$ of real rank zero hypothesis or the $AF$ hypothesis is easily dealt with by this method.

(iv) In some other cases the separability hypothesis on $I$ can be removed by reducing to the separable case. For example, this will work if $I$ is the tensor product of an elementary $C^*$--algebra with $C(T)$, $T$ compact, Hausdorff, and finite dimensional. But $I$ don't know how to remove the separability hypothesis in general.

It is probably premature to define bootstrap categories, so the next corollary should be regarded as just an illustration. In particular the category $\mathcal C$ could already be enlarged, at the cost of having separate categories for real and stable rank, by using parts (iii) and (iv) of Remark 3.13.

\medskip
\noindent
{\bf Corollary\ 3.14.\ }\textit{Let $\mathcal C$ be the smallest class of $C^*$--algebras containing all those satisfying the hypotheses on $I$ in either 3.11 on 3.12 and such that:}
\item{(i)}\ \textit{If $I$ is an ideal of $B$ such that both $I$ and $B/I$ are in $\mathcal C$, then $B$ is in $\mathcal C$,}
\item{(ii)\ }\textit{If prim$(B)=\bigcup^\infty_{n=1} F_n$ where each $F_n$ is closed, and if each $B(F_n)$ is in $\mathcal C$, then $B$ is in $\mathcal C$, and}
\item{(iii)\ }\textit{If $B=(\bigcup J_\la)^-$ where $\{J_\la\}$ is an upward directed family of ideals, and if each $B/J_\la^\perp$ is in $\mathcal C$, then $B$ is in $\mathcal C$.}

\medskip
\noindent
\textit{Then if $I$ is an ideal of a $C^*$--algebra $A$ and $I$ is in $\mathcal C$, we have $\rt(A)\leq\rma(\rt(I),\rt(A/I),2)$ and $\rR(A)=\rma(\rR(I),\rR(A/I))$.}

\begin{proof}The validity of (ii) follows from (CST) as in the first part of the proof of 3.12. And the validity of (iii) follows from Sheu's Technical Proposition, [Sh, 3.15], and its real rank counterpart, [O1, Lemma 1.9]. Note that (iii) is a special case of (ii) when $B$ is separable, since then $\{J_\la\}$ may be assumed countable.
\end{proof}

Since much generalization of the results of Dixmier and Douady [DD] has been 
done, one hopes that Theorem 3.12 can be generalized.  The following uses the 
corona factorization property, a concept introduced by Kucerovsky and Ng, 
cf. [DN, Definition 2.1], to abstract the key part of the proof of 3.12.

\medskip
\noindent
{\bf Corollary\ 3.15.\ }\textit{Assume that $I$ is a separable stable ideal 
of the $C^*$--algebra $A$ and that $I$ has the corona factorization 
property and has an approximate identity consisting of projections.  Then}
\item{(i)}\ $\rt(A) \le \rma(\rt(A/I),\;2),\;$\textit{and}
\item{(ii)\ }$\rR(A) \le \rma(\rR(A/I),\;1).$

\begin{proof}(ii) The proof proceeds like that of 3.12 through the 
construction of the projections $p$ (in $I$) and $q$ (in $\mathbb M_n(A)$), 
but $p$ is no longer full.  Let $B=(\bon_n-q)\mathbb M_n(I)(\bon_n-q)$.  
The key point is to prove that $B$ is stable, and we first prove that 
$\bon_n-q$ is full in $\mathbb M_n(A)$.  In fact if $J$ is a proper ideal 
of $A$ such that $\bon_n-q\in \mathbb M_n(J)$, then Lemma 2.6 can be 
applied in $A/J =\lambda(A)$ to obtain a contradiction.  Note that 
\[\uy(\bon-p)\uy^*\ge\delta q \Rightarrow\lambda(\uy\,\uy^*)\ge\lambda(\uy(\bon-p)\uy^*)\ge\delta\cdot\bon_{\mathbb M_n(A/J)}.\]
Then it is easy to deduce from [KN, Definition 2.1] that $B$ is stable.  It 
then follows from [Br1, Theorem 3.1] or [K, Theorem 2] (cf. [Br2, Theorem 4.23 and page 963]) that there exists a subprojection $r$ of $\rho_n(\bon_n-q)$ such that $r=\sum_1^{\infty} r_m$, where the $r_m$'s and the sum are as in the proof of 
3.12.  The rest of the proof is just like that of 3.12.

(i) A proof can be given which is like that of (ii) with two exceptions:
\item{1.}\ The substitute for Lemma 2.6 is provided by [Ri].  First, we 
know \textit{a priori} that $\rt(A) <\infty$.  And thus [Ri, Proposition 6.5] 
implies that no non-trivial quotient of $A$ can have an $n$--tuple $\underline w$ with $n>1$ and $\underline w\,\underline w^*$ invertible.
\item{2.\ }The $\wtk$ term can be omitted.
\end{proof}

\medskip
\noindent
{\bf 3.16.\ Multiplier Algebras.}

\noindent{\bf (i) Example.}  There is a separable subhomogeneous $C^*$--algebra $A$ such that $\rt(M(A))>\rt(A)$ and $\rR(M(A))>\rR(A)$. It is also true that prim$(A)=\bigcup^\infty_{n=1}\,F_n$, where each $F_n$ is closed and each $A(F_n)$ is unital. Thus this example shows that cases (i$'$) and (ii$'$) of [BP2, Theorem 2.12] cannot be extended to higher ranks, justifying a claim made in [BP2, Remark 2.13(iii)]. Let $X$ be a ball of dimension $d\geq 4$ and $n$ a positive integer such that $n\geq(d+3)/2$. Thus
\[ \rt(C(X)\otimes \mathbb M_n)=\rt(C(X)\otimes \mathbb M_{n-1})=2,\; \mathrm{and}\]
\[ \rR(C(X)\otimes \mathbb M_n)=\rR(C(X)\otimes \mathbb M_{n-1})=1.\]
 Let $B_1=\mathbb M_n(C(X))$ and $B_0=\{(f_{ij})\in B_1:f_{in}=f_{nj}=0\}$. Thus $B_0\cong \mathbb M_{n-1}(C(X))$. Finally, let $A=A_d=$
\[\{(a_m)^\infty_{m=1}:a_m\in B_1, \forall m,\;\mathrm{and}\; (a_m)\; \mathrm{converges\; to\; an\; element\; of\;} B_0\}.\]
 Then $\rp(A)=\cup_{1\leq m\leq\infty}F_m$, where $F_m= X$ for all $m$, $A(F_m)\cong B_1$ for $m<\infty$, and $A(F_\infty)\cong B_0$. In particular, by (CST), $\rt(A)=2$ and $\rR(A)=1$.

Represent elements of $B_1$ as $2\times 2$ block matrices, $\begin{pmatrix}a & b\\c & d\end{pmatrix}$, where $a$ is $(n-1)\times(n-1)$ and $d$ is $1\times 1$. Then $M(A)$ can be identified with the set of bounded sequences $\left(\begin{pmatrix}a_m & b_m\\c_m & d_m\end{pmatrix}\right)$ such that $(a_m)$ is convergent and $(b_m), (c_m)$ converge to 0. Let $p$ be the constant sequence $\left(\begin{pmatrix}{\bf{1_{n-1}}} & 0\\0 & 0\end{pmatrix}\right)$, and note that $p\in A$. Thus $M(A)/A\cong(\bon-p)M(A)(\bon-p)/(\bon-p)A(\bon-p)$. Since $(\bon-p)M(A)(\bon-p)$ is the $\ell^\infty$--direct sum (or direct product) of countably infinitely many copies of $C(X)$ and $(\bon-p)A(\bon-p)$ is the $c_0$--direct sum, it is easily seen that $\rt(M(A)/A)\geq\rt(C(X))>2$ and $\rR(M(A)/A)\geq\rR(C(X))>1$. The only technical point involved in verifying this last assertion is to note that if $\ux=(\ux_m)$ is a tuple in the $\ell^\infty$--direct sum whose image in the quotient is unimodular, then $\ux_m$ is unimodular for all but finitely many values of $m$.  (In fact, $\rt(M(A)/A)=\rt(C(X))$ and $\rR(M(A)/A)=\rR(C(X))$,
cf. (ii) below.)

Finally let $C$ be the $c_0$--direct sum of the algebras $A_d$ constructed above for all values of the dimension $d$. Thus $C$ is still separable and has only finite dimensional irreducible representations, but $C$ is no longer subhomogeneous, and $\topdim(C)=\infty$. Then $\rt(M(C)=\rR(M(C))=\infty$. But $\rp(C)=\bigcup^\infty_{m=1} F_m$, where each $F_m$ is closed, each $C(F_m)$ is unital, $\rt(C(F_m))=2$, and $\rR(C(F_m))=1$. In particular, $\rt(C)=2$ and $\rR(C)=1$.

\medskip
\noindent
{\bf (ii)\ Remark.\ }On the other hand, if $A$ is separable and subhomogeneous, then $\topdim(M(A))=\topdim(A)$. (Techniques for reduction to the separable 
case allow the separability hypothesis to be weakened to $\sigma$--unitality, 
but the argument is longer than most of this type and will be omitted.)  I think this result should be essentially known, possibly folklore, but haven't been able to find a reference. The first step is to prove the following, which I learned from conversations with M. Dupre in the 1970's:

If $A$ is separable (or just $\sigma$--unital) and $n$--homogeneous, and if top dim$(A)<\infty$, then $M(A)$ is $n$--homogeneous and $\rp(M(A))=\beta(\rp(A))$, the Stone--\v{C}ech compactification.

Let $X=\rp(A)$, so that $X$ is $\sigma$--compact, locally compact, Hausdorff, and finite dimensional, and $A$ is given by a locally trivial $\mathbb M_n$--bundle on $X$. This bundle is necessarily of finite type, in the sense that $X$ can be covered by finitely many (actually, $1+\rd(X)$) open sets over each of which the bundle is trivial ([Hu, 3.5.4]).  The facts that the bundle is of finite type and that $\mathrm{Aut}(\mathbb M_n)$ is compact, combined with standard techniques relating to Stone--\v{C}ech compactifications of normal spaces, allow one to extend the bundle to $\beta(X)$. Once one has a bundle over $\beta(X)$, similar techniques show that $M(A)$ consists of the bounded sections of this bundle.

Now since $X$ is normal, $\topdim(M(A))=\rd(\beta(X))=\rd(X)$; and since $X$ is $\sigma$--compact, $\rd(X)=\mathrm{loc}\;\rd(X)=\topdim(A)$. This covers the case where $A$ is homogeneous, and the general case is proved by induction on $n$, where $n$ is the maximum dimension of an irreducible representation. There is an ideal $I$ which is $n$--homogeneous such that all irreducibles of $A/I$ have dimension less than $n$. By the non--commutative Tietze extension theorem, whose separable case is [P, 3.12.10], the natural map $\opi:M(A)\to M(A/I)$ is surjective. The kernel, $M(A,I)$, of $\opi$ is isomorphic to a hereditary $C^*$--subalgebra of $M(I)$. Since $\rp(M(A,I))$ is an open subset of $\rp(M(I))$, $\topdim(M(A,I))\leq\topdim(M(I))$ (but $\rra(M(A,I))$ may be much bigger than $\rra(M(I))$. So the induction goes through easily. 

\medskip
\noindent
{\bf (iii)\ Questions.\ }Both (i) and (ii) above relate to the desire for non--commutative analogues of the theorem that dim$(\beta(X))=\rd(X)$ for a normal topological space $X$. It is natural to draw the conclusion that, in the higher rank situation, one should focus on $\topdim{}$ rather than real or stable rank. However, [BP2, Corollary 3.10] includes a positive result about $\rt(M(A))$ under special hypotheses. And more importantly, except in the zero--dimensional case, $\topdim(A)$ is defined only when $\rp(A)$ is almost Hausdorff, so that $\topdim(M(A))$ will typically be undefined. (Of course $\mathrm{d_r}(A)=\topdim(A)$, by [W], when $A$ is subhomogeneous, where $\mathrm{d_r}$ is the decomposition rank of Kirchberg and Winter. But $\mathrm{d_r}(M(A))$ will also typically be undefined.) Nevertheless, there are at least two questions on this topic which seem worthy of investigation. Although positive answers would be pleasing, these questions are not conjectures.

\medskip
\noindent
1. If $A$ is a separable $C^*$--algebra all of whose irreducible representations are finite dimensional, is it necessarily true that $\rR(M(A))\leq \topdim(A)$ and $\rt(M(A))\leq 1+\topdim(A)/2$?
\medskip

\noindent
2. Can one prove $\rR(M(A))\leq 1$, or even $\rR(M(A))<\infty$, for $A$ in a reasonably large class of stable $C^*$--algebras?

\medskip
\noindent
Of course, it follows from [Ri, Proposition 6.5] that $\rt(M(A))=\infty$ when $A$ is stable, but the real rank case seems unclear.

\section{A Theorem of P. Green.}

In the theorem below condition (iv) is Dixmier's topological characterization of GCT algebras, and (ii) and (iii) are just intermediate conditions, (iii) being related to Dixmier's concept of Hausdorff point, cf. [D1]. Thus the equivalence of (i)--(iv) is valid without separability. Also one direction of the corollary, that an $FD$--like decomposition implies GCT, is valid without separability and is essentially due to Dixmier. Conditions (v), (vi), and(vii) are new topological characterizations of GCT due to Green. I have made some changes from the presentation of the theorem provided in [G], the only significant one being that the proof given is less topological than the original. In fact the equivalence of the conditions (ii)--(vii) can be proved topologically. Although [G] asserts that all the topological arguments are ``easy'', in one case the best topological argument I could find was not quite easy (though not so terribly hard). Finally, a cover $\{U_i\}$ of a space $X$ is called \textit{point--finite} if no point of $X$ is contained in infinitely many $U_i$'s.

\medskip
\noindent
{\bf Theorem\ 4.1 (Green [G]).} \textit{ If $A$ is a separable $CCR$ $C^*$--algebra, then the following are equivalent:}
\item{(i)\ }\textit{$A$ has generalized continuous trace.}
\item{(ii)\ }\textit{Every non--empty closed subset $F$ of $\rp(A)$ has a non--empty relatively open subset $G$ such that $G$ is Hausdorff and every compact subset of $G$ is closed in $\rp(A)$.}
\item{(iii)\ }\textit{Every non--empty closed subset $F$ of $\rp(A)$ has a non--empty relatively open subset $G$ such that if $x\in G, y\in F$, and $x\neq y$, then $x$ and $y$ have disjoint neighborhoods relative to $F$.}
\item{(iv)\ }\textit{Every non--empty closed subset $F$ of $\rp(A)$ has a non--empty relatively open subset $G$ such that each point of $G$ has a (relative) neighborhood base consisting of sets closed in $\rp(A)$.}
\item{(v)\ }\textit{One can write $\rp(A)=\bigcup^\infty_1 F_n$, where $\{F_n\}$ is a countable family of closed compact sets.}
\item{(vi)\ }\textit{The space $\rp(A)$ is metacompact; i.e., every open cover has an open point--finite refinement.}
\item{(vii)\ }\textit{There is a point--finite open cover $\{U_i\}$ of $\rp(A)$ such that each $U_i$ is contained in a compact subset of $\rp(A)$.}
\item{(viii)\ }\textit{$A$ is stably isomorphic to a $C^*$--algebra with only finite dimensional irreducible representations.}

\begin{proof} (i)$\Rightarrow$(ii): Let $\{J_\alpha:0\leq\alpha\leq\beta\}$ be the continuous trace composition series for $A$, discussed above in 2.4, and let $V_\alpha=\rp(J_\alpha)\subset\rp(A)$. Let $\gamma$ be the smallest index such that $V_\gamma\cap F\neq\emptyset$. Then $\gamma$ cannot be a limit ordinal. Let $G=V_\gamma\cap F=(V_\gamma\setminus V_{\gamma-1})\cap F$. As noted in 2.4, $V_\gamma\setminus V_{\gamma-1}$ has the properties required for $G$.

(ii)$\Rightarrow$(iii): Use the same $G$ produced by (ii). If $y\in G$, the condition follows since $G$ is Hausdorff. If $y\notin G$, the condition follows since $G$ is locally compact.

(iii)$\Rightarrow$(iv): Use the same $G$ produced by (iii), which is necessarily locally compact and Hausdorff. The usual proof that compact subsets of a Hausdorff space are closed now shows that compact subsets of $G$ are closed in $F$, hence globally closed.

(iv)$\Rightarrow$(v): We construct a strictly increasing family $\{V_\alpha:0\leq\alpha\leq\beta\}$ of open sets such that $V_0=\emptyset, V_\beta=\rp(A), V_\la=V_{\alpha<\la} V_\alpha$ if $\la$ is a limit ordinal, and $V_{\alpha+1}\setminus V_\alpha$ has the property specified for $G$ in (iv) relative to $\rp(A)\setminus V_\alpha$ for each $\alpha<\beta$. Since $\rp(A)$ is second countable, the ordinal number $\beta$ is countable. Since $V_{\alpha+1}\setminus V_\alpha$ is second countable and locally quasi--compact for $\alpha<\beta$, each $V_{\alpha+1}\setminus V_\alpha$ is a countable union of closed compact sets.

(v)$\Rightarrow$(vi): We may assume the given family $\{F_n\}$ is increasing. If $\{U_i\}$ is an open cover, choose for each $n$ a finite subcover, $\{V_{nj}:1\leq j\leq m_n\}$, of $F_n$. If $W_{nj}=V_{nj}\setminus F_{n-1}$, then $\{W_{nj}\}$ is an open point--finite refinement of $\{U_i\}$.

(vi)$\Rightarrow$(vii): This is obvious since $\rp(A)$ is second countable and locally quasi--compact.

(vii)$\Rightarrow$(viii): The point here is that the open cover $\{U_i\}$ provided by (vii) makes it possible to do a better version of the argument for [Br1, 2.11 a]. Since $\rp(A)$ is second countable, we may assume $\{U_i\}$ is countable. Let $U_i\subset K_i$, $K_i$ compact. For each $P$ in $K_i$ choose $e_P\in K(A)_+$, where $K(A)$ is the Pedersen ideal, such that ${e_P}\notin P$. If $V_P=\{Q\in\rp(A):e_P\notin Q\}$, then the $V_P$'s form an open cover of $K_i$, and there is a finite subcover $\{V_{P_j}:1\leq j\leq m_i\}$. Let $f_i=\sum^{m_i}_1\,e_{P_j}$. Thus $f_i\in K(A)_+$ and $f_i\not\in P$ for $P$ in $K_i$ or, \textit{a fortiori}, 
for $P$ in $U_i$. Now let $I_i$ be the ideal $A(U_i)$, and let $g_i$ be a strictly positive element of $I_i$. Then $f_ig_if_i\in I_i\cap K(A)_+$, and $f_ig_if_i$ generates $I_i$ as an ideal. (If $\pi$ is an irreducible representation such that $\mathrm{ker}\; \pi\in U_i=\rp(I_i)$, then $\pi(f_i)\neq 0$ and $\pi(g_i)$ is a positive operator with trivial nullspace.) Now let $h=\Sigma_i \,\epsilon_if_ig_if_i$, where the $\epsilon_i$'s are positive numbers such that $\Sigma\, \epsilon_i\|f_ig_if_i\|<\infty$. Then $\pi(h)$ is a non--zero finite rank operator for each irreducible representation $\pi$ of $A$, since $\{U_i\}$ is a point--finite cover. If follows that $B=(h A h)^-$ is a full hereditary $C^*$--subalgebra of $A$ all of whose irreducible representations are finite dimensional. By [Br1] $B$ is stably isomorphic to $A$.

(viii)$\Rightarrow$(i): This is essentially due to Dixmier, but I sketch a proof. If $A\otimes \bK\cong B\otimes \bK$, where $B$ has only finite dimensional irreducibles, then a Baire category argument produces a composition series $\{I_\alpha:0\leq\alpha\leq\beta\}$ for $B$ such that for each $\alpha<\beta$, $I_{\alpha+1}/ I_\alpha$ is $n_\alpha$--homogeneous for some natural number $n_\alpha$. If $\{I'_\alpha\}$ is the corresponding composition series for $A$, then $I'_{\alpha+1}/ I'_\alpha\subset J(A/I'_\alpha)$.
\end{proof}

\medskip
\noindent
{\bf Corollary\ 4.2.\ }\textit{If $A$ is a separable $CCR$ $C^*$--algebra, then $A$ has generalized continuous trace if and only if $\rp(A)$ has an $FD$--like decomposition.}

\begin{proof} If $\{H_n\}$ is an $FD$--like decomposition of $\rp(A)$ then each $H_n$ is the union of countably many compact sets, each of which is necessarily closed in $\rp(A)$. Thus condition (v) is satisfied. The converse follows directly from (i)$\Rightarrow$(viii).
\end{proof}

\bigskip

\noindent Department of Mathematics

\noindent Purdue University

\noindent West Lafayette

\noindent Indiana 47907, USA

\noindent lgb@math.purdue.edu

\end{document}